\documentclass[12pt,letterpaper]{article}

\UseRawInputEncoding
\usepackage{bbm}
\usepackage{chngcntr}
\usepackage{multirow}

\usepackage[colorlinks=true,urlcolor=blue,citecolor=blue,linkcolor=blue,bookmarks=true,bookmarksopen=false]{hyperref}

\usepackage{algorithmic}
\usepackage[linesnumbered,ruled,vlined]{algorithm2e}

\usepackage[left=1in,top=1in,right=1in,bottom=1in]{geometry}

\newcommand{\exclude}[1]{}

\renewcommand{\Omega}{\varOmega}

\newcounter{commentcounter}
\long\def\symbolfootnote[#1]#2{\begingroup%
\def\thefootnote{\fnsymbol{footnote}}\footnote[#1]{#2}\endgroup}

\usepackage{sibarticle}

\usepackage{amsmath}
\usepackage{kbordermatrix}
\usepackage{mathtools}

\usepackage{tikz}
\usetikzlibrary{arrows,shapes,snakes,automata,backgrounds,petri}
\usetikzlibrary{graphs}

\newcommand{\ignore}[1]{}

\usepackage{color}
\usepackage{bm}
\usepackage{multirow}
\usepackage{booktabs}

\usepackage{rotating}
\usepackage[symbol]{footmisc}
\usepackage{graphicx}
\usepackage{subfigure}
\usepackage{epstopdf}

\allowdisplaybreaks

\title{A Note on Valid Inequalities for PageRank Optimization with Edge Selection Constraints}
\ShortTitle{A Note on Valid Inequalities for PageRank Optimization with Edge Selection Constraints} \ShortAuthors{Authors}
\NumberOfAuthors{1} \FirstAuthor{Shang-Ru Yang, Yung-Han Liao, Chih-Ching Chien, Hao-Hsiang Wu\footnotemark[1]}
\FirstAuthorAddress{ Department of Management Science, National Yang Ming Chiao Tung University, Hsinchu, Taiwan\\ 
	{\tt \{brian895.mg11, lyh001115.mg12, ccchien11.c, hhwu2\}@nycu.edu.tw	}}

\SecondAuthor{Authors}
\SecondAuthorAddress{
	TBD\\ {\tt
		TBD}
}
\ShortAuthors{Yang, Liao, Chien, and Wu}
\exclude{

}

\renewcommand{\baselinestretch}{1.25}


\begin{document}

\maketitle \centerline{\today}
\footnotetext[1]{Corresponding author (hhwu2@nycu.edu.tw)} 
\begin{abstract}
Cs\'{a}ji, Jungers, and Blondel prove that while a PageRank optimization problem with edge selection constraints is NP-hard, it can be solved optimally in polynomial time for the unconstrained case. This theoretical result is accompanied by several observations, which we leverage to develop valid inequalities in polynomial time for this class of NP-hard problems.
\end{abstract}

\renewcommand{\baselinestretch}{1.5}

\section{Introduction}\label{problem}

This note considers the NP-hard PageRank optimization problem with edge selection constraints from Cs\'{a}ji, Jungers, and Blondel \cite{PageRankSSP1,PageRankSSP}. Let $G = (V,v,E)$ be a directed graph with a set of nodes $V = \{1,\dots,n\}$ for $n \in \mathbb N$, a target node $v \in V$, and a set of directed edges $E$. Suppose that there exists a set of hidden (fragile) edges $Z$ in $G$ such that $Z \cap E = \emptyset$, making $Z$ and $E$ mutually exclusive. Given a graph $G=(V,v,E)$, and a set of hidden edges $Z$, our goal is to identify a subset $Y \subseteq Z$ under specific constraints such that $v$ attains an optimal PageRank value of \cite{PageRRank1998} in the graph $G=(V, v, E \cup Y)$.  Previous studies, such as those by \cite{PLink2006, Kerchove2008}, explore the optimal structure for outgoing edges from certain nodes in $V$. Similarly, another study \cite{Olsen2014} explores the theoretical difficulty of selecting a fixed number of incoming edges to the target node $v$ to optimize its PageRank. The result demonstrates that, unless P = NP, there is no polynomial-time approximation method to maximize PageRank of $v$ by selecting the fixed number of incoming edges. 

The PageRank value of a node $v$ equals the inverse of the expected {\it first return} time to $v$, i.e., the average steps for a random walk starting at $v$ to return to $v$ in $G$ \cite{Aldous2002, Sheldon2007}. Cs\'{a}ji et al. study maximizing the PageRank value of a node $v$ by equivalently minimizing the expected first return time to $v$ through the selection of edges in $Z$ \cite{PageRankSSP1,PageRankSSP}. Formally, we define $\mathcal{FR}({Y})$ as the expected first return time to $v$ after selecting (adding) the set of edges $Y \subseteq Z$ to the graph $G$. For each directed edge $(i,j) \in Z$, we define $y_{(i,j)} \in \mathbb B$ as a binary variable, where $y_{(i,j)} = 1$ denotes the selection of the edge $(i,j) \in Z$; otherwise, $y_{(i,j)} = 0$. Specifically, an incumbent selection $\bar Y \subseteq Z$ is equivalent to the support $\bar Y  = \{(i,j) \in Z: \bar y_{(i,j)} = 1 \}$. For the notation $\bar Y \subseteq Z$ and vector $\bar {\bold{y}} \in \mathbb B^{|Z|}$, we make a slight adjustment for convenience. We name the notation $\bar {\bold{y}} \in \mathbb B^{|Z|}$ and its support $\bar Y$, and name to the corresponding function evaluations $\mathcal{FR}(\bar {\bold{y}})$ for  $\bar {\bold{y}} \in \mathbb B^{|Z|}$ and $\mathcal{FR}(\bar Y)$ for the corresponding support $\bar Y \subseteq {Z}$,  interchangeably. We follow the concept considered by the study of Cs\'{a}ji et al., which minimizes $\mathcal{FR}({Y})$ based on the selected edges $Y \subseteq Z$. Let $\mathcal{Y}$ be a set of constraints on the associated $|Z|$-dimensional binary decision vector $\bold{y} \in \mathbb B^{|Z|}$.  Given a graph $G = (V,v,E)$ and a set $Z$, the general PageRank optimization problem with edge selection constraints from Cs\'{a}ji et al. is defined as
\begin{equation}\label{SSP_E}
	\min_{\bold{y} \in \mathcal{Y} \cap \mathbb B^{|Z|}} \mathcal{FR}(\bold{y}).
\end{equation} 
We note that Problem \eqref{SSP_E} presents the potential to explore various structures in mathematical programming. For instance, if the evaluation of an objective function relies on multiple sampling scenarios based on the sample average approximation method of \cite{KSH02}, stochastic programming approaches such as those in \cite{BL97,Higle1991,SMIP2017,SDR09} can be considered. Similarly, if probabilistic constraints are in $\mathcal Y$, the chance-constrained optimization methods \cite{AhmedCCP2017,SimgeCCP2022,Luedtke2008} could be explored. Based on general methods of integer programming \cite{MINLP2013,Conforti2014,NW88,Schrijver98}, Problem \eqref{SSP_E} can be solved using an exact cutting plane method based on a relaxed master problem (RMP)
\begin{equation}\label{RMP}
	\min \{\theta: (\theta, \bold{y}) \in \mathcal{C}, \bold{y} \in \mathcal{Y}\cap \mathbb B^{|Z|}, \theta \in \mathbb R_+ \},
\end{equation}
where $\theta \in  \mathbb R_+$ is a variable that evaluates the value of the $\mathcal{FR}$ function and $\mathcal{C}$ is a set of valid inequalities iteratively updated to refine the RMP \eqref{RMP}. At each iteration, given an incumbent solution $\bold{\bar y} \in \mathcal{Y} \cap \mathbb B^{|Z|}$, the method generates a valid inequality associated with this incumbent $\bold{\bar y}$ and adds it to $\mathcal{C}$. The process continues until the optimality guarantee is achieved, ensuring convergence to the optimal solution.

This note considers valid inequalities for the set $\mathcal{C}$ of the RMP \eqref{RMP}. We note that Cs\'{a}ji et al. prove that the unconstrained case $\mathcal Y = \emptyset$ of Problem \eqref{SSP_E} can be solved optimally in polynomial time \cite{PageRankSSP1,PageRankSSP}. Later, a similar conclusion is independently reached by \cite{Fercoq2013}, further confirming the validity of this result. We summarize the following observations from the formal results of \cite{PageRankSSP}, which can help us develop valid inequalities.
\begin{observation} \label{OB}
	From the results of \cite{PageRankSSP}, we obtain the following observations.
	\begin{itemize}		
		\item[$(i)$] The problem $\min_{\bold{y} \in  \mathbb B^{|Z|}} \mathcal{FR}(\bold{y})$ can be solved in polynomial time.
		\item[$(ii)$]Given two exclusive subsets $\bar S \subseteq Z$ and $\bar N \subseteq Z$, the  function 
		\begin{equation}\label{gamma_f}
			\gamma (\bar S, \bar N) = \min \{ \mathcal{FR}(\bold{y}): y_{(i,j)} = 1 \quad \forall (i,j) \in \bar S,\ y_{(i,j)} = 0 \quad \forall (i,j) \in \bar N,\ \bold{y} \in  \mathbb B^{|Z|} \}
		\end{equation}
		can be solved in polynomial time.   
		\item[$(iii)$] Given an incumbent solution $\bold{\bar y} \in \mathcal{Y} \cap \mathbb B^{|Z|}$, the function $\mathcal{FR}(\bold{\bar y})$ can be solved in polynomial time, where $\mathcal{FR}(\bold{\bar y}) \ge 0$. 
	\end{itemize}
\end{observation}
We take the $L$-shaped inequality from \cite{LL93} as an example of the baseline for the set $\mathcal{C}$. Let $L \in \mathbb R$ be a small enough constant that makes the $L$-shaped inequality valid. Given an incumbent $\bar{\bold{y}} \in \mathcal{Y} \cap \mathbb B^{|Z|}$, the $L$-shaped inequality is defined as 
\begin{equation*}
\theta \ge \mathcal{FR}(\bar {\bold{y}}) + \sum_{(i,j)\in \bar Y}  (L-\mathcal{FR}(\bar {\bold{y}}))(1-y_{(i,j)})+\sum_{(i,j)\in Z\setminus \bar Y} (L-\mathcal{FR}(\bar {\bold{y}})) y_{(i,j)}.
\end{equation*}	
Here, since $\mathcal{FR}(\bold{\bar y}) \ge 0$ for any $\bold{\bar y} \in \mathcal{Y}\cap \mathbb B^{|Z|}$, constant $L = 0$  is trivially and sufficiently small to provide a valid inequality for the RMP \eqref{RMP}; however, through Observation \ref{OB}, we can compute a stronger value $\min_{\bold{y} \in \mathbb B^{|Z|}} \mathcal{FR}(\bold{y}) \ge 0$ polynomially  for $L$ to have a stronger $L$-shaped inequality
\begin{equation}\label{L_cut}
	\theta \ge \mathcal{FR}(\bar {\bold{y}}) + \sum_{(i,j)\in \bar Y}  (	\min_{\bold{y} \in \mathbb B^{|Z|}} \mathcal{FR}(\bold{y})-\mathcal{FR}(\bar {\bold{y}}))(1-y_{(i,j)})+\sum_{(i,j)\in Z\setminus \bar Y} (\min_{\bold{y} \in \mathbb B^{|Z|}} \mathcal{FR}(\bold{y})-\mathcal{FR}(\bar {\bold{y}})) y_{(i,j)}.	
\end{equation}
We note that, without  Observation  \ref{OB}, even the $L$-shaped inequality \eqref{L_cut} cannot be guaranteed to be generated in polynomial time. Starting with the $L$-shaped inequality \eqref{L_cut}, and based on Observation  \ref{OB}, we develop stronger valid inequalities in the following sections.

\section{A Valid Inequality} \label{method}

Given an incumbent $\bar Y$, we can utilize the function \eqref{gamma_f} of Observation \ref{OB} to calculate the best case of either adding an edge $(i,j)$ from $Z\setminus \bar Y$ or not adding an edge $(i,j)$ from $\bar Y$. Note that base on the function \eqref{gamma_f}, together with the condition $\mathcal{FR}(\bold{\bar y}) \ge 0$ for any $\bold{\bar y} \in \mathcal{Y}\cap \mathbb B^{|Z|}$, we develop another class of valid inequalities for the set $\mathcal C$ of the RMP \eqref{RMP}.
\begin{theorem}\label{prop:New_valid}
	Given an incumbent $\bar{\bold{y}} \in \mathcal{Y} \cap \mathbb B^{|Z|}$, the inequality 
\begin{align}\label{New_cut}
		\theta \ge \mathcal{FR}(\bar {\bold{y}}) & + \sum_{(i,j)\in \bar Y}  \min\{0, (\gamma (\emptyset, \{(i,j)\})-\mathcal{FR}(\bar {\bold{y}})) \} (1-y_{(i,j)})\nonumber \\ 
		&+\sum_{(i,j)\in Z\setminus \bar Y} \min\{0, (\gamma ( \{(i,j)\}, \emptyset)-\mathcal{FR}(\bar {\bold{y}})) \} y_{(i,j)}, 
\end{align} 	
	is valid for the RMP \eqref{RMP}.
\end{theorem}    
\begin{proof}
	Consider a feasible point $(\hat {\bold{y}}, \hat \theta)$ to the RMP \eqref{RMP}. We show that $(\hat {\bold{y}}, \hat \theta)$ satisfies the inequality \eqref{New_cut} for any $\bar{\bold{y}} \in \mathcal{Y} \cap \mathbb B^{|Z|}$ with the associated $\hat Y \neq \bar Y$. Because of $\hat Y \neq \bar Y$, there exists two cases: at least one edge of $\bar Y$ that is not in $\hat Y$, or at least one edge of $\hat Y$ that is not in $ \bar Y$ for the inequality \eqref{New_cut}. 
	\begin{itemize}		
	\item[Case 1:] For the case where at least one edge of $\bar Y$ is not in $\hat Y$, there are two possible results in the term $(\gamma (\emptyset, \{(i,j)\})-\mathcal{FR}(\bar {\bold{y}}))(1-\hat y_{(i,j)})$ when $\hat y_{(i,j)}=0$ for an edge $(i,j) \in \bar Y \setminus \hat Y$. We let $(i,j)_0'$ and $(i,j)_0''$ denote two  edges in $\bar Y \setminus \hat Y$ for the two results, where the associated variables $\hat y_{(i,j)_0'} = \hat y_{(i,j)_0''} = 0$ in $\hat {\bold{y}}$ with the conditions $(\gamma (\emptyset, \{(i,j)_0'\})-\mathcal{FR}(\bar {\bold{y}}))(1-\hat y_{(i,j)_0'}) < 0$ and $(\gamma (\emptyset, \{(i,j)_0'' \})-\mathcal{FR}(\bar {\bold{y}}))(1-\hat y_{(i,j)_0''}) \ge 0$. Here, since $\hat {\bold{y}}$ must includes at least one $\hat y_{(i,j)_0'} =0 $ or $\hat y_{(i,j)_0''} = 0$, we have $\gamma (\emptyset, \{(i,j)_0''\}) = \min \{  \mathcal{FR}(\bold{y}): y_{(i,j)_0''} = 0,  \bold{y} \in \mathbb B^{|Z|}\} \le \mathcal{FR}(\hat {\bold{y}})$ and $\gamma (\emptyset, \{(i,j)_0'\}) = \min \{  \mathcal{FR}(\bold{y}): y_{(i,j)_0'} = 0,  \bold{y} \in \mathbb B^{|Z|}\} \le \mathcal{FR}(\hat {\bold{y}})$. In other words, both $\gamma (\emptyset, \{(i,j)_0''\})$ and $\gamma (\emptyset, \{(i,j)_0'\})$ represent the relaxation of $\mathcal{FR}(\hat {\bold{y}})$.
	\item[Case 2:] Similar to the discussion in the previous case, we have another situation that at least one edge of $\hat Y$ is not in $\bar Y$. To put it briefly, let $(i,j)_1'$ and $(i,j)_1''$ be possible edges in $\hat Y \setminus \bar Y$, where the associated variables $\hat y_{(i,j)_1'} = \hat y_{(i,j)_1''} = 1$ of $\hat {\bold{y}}$ with $(\gamma (\{(i,j)_1'\},\emptyset )-\mathcal{FR}(\bar {\bold{y}}))\hat y_{(i,j)_1'} < 0$ and $(\gamma (\{(i,j)_1''\},\emptyset )-\mathcal{FR}(\bar {\bold{y}}))\hat y_{(i,j)_1''} \ge 0$. Both $\gamma (\{(i,j)_1''\},\emptyset )=\min \{  \mathcal{FR}(\bold{y}): y_{(i,j)_1''} = 1,  \bold{y} \in \mathbb B^{|Z|}\} \le \mathcal{FR}(\hat {\bold{y}})$ and $\gamma (\{(i,j)_1'\},\emptyset )=\min \{  \mathcal{FR}(\bold{y}): y_{(i,j)_1'} = 1,  \bold{y} \in \mathbb B^{|Z|}\} \le \mathcal{FR}(\hat {\bold{y}})$ are the relaxation of $\mathcal{FR}(\hat {\bold{y}})$. 
\end{itemize}
We demonstrate the validity of the following inequalities for a point $(\hat {\bold{y}}, \hat \theta)$ based on the discussion in Cases 1 and 2.
	\begin{align}
		\hat \theta & \ge \mathcal{FR}(\hat {\bold{y}})\nonumber\\
		& \ge  \mathcal{FR}(\bar {\bold{y}}) +  \min\{0, (\gamma ( \emptyset,\{(i,j)_0''\}) -\mathcal{FR}(\bar {\bold{y}})) \} (1 - \hat y_{(i,j)_0''})+  \min\{0, (\gamma (\{(i,j)_1''\},\emptyset )-\mathcal{FR}(\bar {\bold{y}})) \} \hat y_{(i,j)_1''}\label{New_valid_1}  \\
		& \ge  \mathcal{FR}(\bar {\bold{y}}) +  (\gamma ( \emptyset,\{(i,j)_0'\})-\mathcal{FR}(\bar {\bold{y}})) (1 - \hat y_{(i,j)_0'}) +  (\gamma (\{(i,j)_1'\},\emptyset )-\mathcal{FR}(\bar {\bold{y}}))  \hat y_{(i,j)_1'} \label{New_valid_2}  \\
		& \ge \mathcal{FR}(\bar {\bold{y}}) + \sum_{(i,j)\in \bar Y}  \min\{0, (\gamma ( \emptyset,\{(i,j)\})-\mathcal{FR}(\bar {\bold{y}})) \} (1-\hat y_{(i,j)}) \nonumber \\
		& \quad \quad \quad \quad \   + \sum_{(i,j)\in Z \setminus \bar Y}  \min\{0, (\gamma (\{(i,j)\}, \emptyset)-\mathcal{FR}(\bar {\bold{y}})) \} \hat y_{(i,j)}. \label{New_valid_3}
	\end{align} 
Based on the relations $(\gamma ( \emptyset,\{(i,j)_0''\}) -\mathcal{FR}(\bar {\bold{y}})) \ge 0$ with $\mathcal{FR}(\hat {\bold{y}}) \ge \gamma (\emptyset, \{(i,j)_0''\})$ and $(\gamma (\{(i,j)_1''\},\emptyset )-\mathcal{FR}(\bar {\bold{y}})) \ge 0$ with $\mathcal{FR}(\hat {\bold{y}}) \ge \gamma (\{(i,j)_1''\},\emptyset )$, we conclude $\mathcal{FR}(\hat {\bold{y}}) \ge  \gamma (\emptyset, \{(i,j)_0''\}) \ge \mathcal{FR}(\bar {\bold{y}})$ and $\mathcal{FR}(\hat {\bold{y}}) \ge  \gamma (\{(i,j)_1''\},\emptyset ) \ge \mathcal{FR}(\bar {\bold{y}})$  for the validity of the inequality \eqref{New_valid_1}, where the inequality \eqref{New_valid_1} can also be regarded as $\hat \theta  \ge \mathcal{FR}(\hat {\bold{y}})\nonumber
	\ge  \mathcal{FR}(\bar {\bold{y}}) + 0 \times (1 - \hat y_{(i,j)_0''}) + 0  \times  y_{(i,j)_1''}$. The inequality \eqref{New_valid_2} follows from the  relations $\mathcal{FR}(\hat {\bold{y}}) \ge \gamma (\emptyset, \{(i,j)_0'\}) = \mathcal{FR}(\bar {\bold{y}}) +  (\gamma (\emptyset, \{(i,j)_0'\})-\mathcal{FR}(\bar {\bold{y}})) (1-y_{(i,j)_0'})$ with $(\gamma (\{(i,j)_1'\},\emptyset )-\mathcal{FR}(\bar {\bold{y}}))  < 0$ and $\mathcal{FR}(\hat {\bold{y}}) \ge \gamma (\{(i,j)_1'\},\emptyset ) = \mathcal{FR}(\bar {\bold{y}}) +  (\gamma (\{(i,j)_1'\},\emptyset )-\mathcal{FR}(\bar {\bold{y}}))y_{(i,j)_1'}$ with $(\gamma (\{(i,j)_1'\},\emptyset )-\mathcal{FR}(\bar {\bold{y}})) < 0$. The inequality \eqref{New_valid_3} follows from $\min\{0, (\gamma ( \emptyset,\{(i,j)\})-\mathcal{FR}(\bar {\bold{y}})) \} \le 0$ for all $(i,j) \in \bar Y \setminus \{(i,j)_0'\}$ and $\min\{0, (\gamma (\{(i,j)\}, \emptyset)-\mathcal{FR}(\bar {\bold{y}})) \} \le 0$ for all $(i,j) \in Z \setminus ( \bar Y \cup \{(i,j)_1'\} )$. This completes the proof.
\end{proof} 
Given an incumbent $\bar{\bold{y}}$, the inequality \eqref{New_cut} discusses the best possible outcome for each individual element within the associated $\bar Y$ and $Z \setminus \bar Y$. We provide the following theoretical result, demonstrating that the inequality \eqref{New_cut} is stronger than the $L$-shaped inequality \eqref{L_cut}.
\begin{proposition}\label{prop:strong}
	Given an incumbent $\bar{\bold{y}} \in \mathcal{Y} \cap \mathbb B^{|Z|}$, the inequality \eqref{New_cut}  is stronger than the inequality \eqref{L_cut}.
\end{proposition}    
\begin{proof}
 Since $\min_{\bold{y} \in \mathbb B^{|Z|}} \mathcal{FR}(\bold{y})$ is a relaxation of both $\min \{  \mathcal{FR}(\bold{y}): y_{(i,j)} = 1,  \bold{y} \in \mathbb B^{|Z|}\} = \gamma ( \{(i,j)\},\emptyset)$ and $ \min \{  \mathcal{FR}(\bold{y}): y_{(i,j)} = 0,  \bold{y} \in \mathbb B^{|Z|}\} = \gamma ( \emptyset,\{(i,j)\})$, we have $\gamma ( \{(i,j)\},\emptyset) \ge \min_{\bold{y} \in \mathbb B^{|Z|}} \mathcal{FR}(\bold{y})$ and $\gamma ( \emptyset,\{(i,j)\}) \ge \min_{\bold{y} \in \mathbb B^{|Z|}} \mathcal{FR}(\bold{y})$ for all $(i,j) \in Z$. Because of $\min_{\bold{y} \in \mathbb B^{|Z|}}\mathcal{FR}(\bold{y}) \le \mathcal{FR}(\bar {\bold{y}})$, we have
	\begin{equation*}
		\min_{\bold{y} \in \mathbb B^{|Z|}}\mathcal{FR}(\bold{y})-\mathcal{FR}(\bar {\bold{y}}) \le \min\{0, (\gamma ( \emptyset,\{(i,j)\})-\mathcal{FR}(\bar {\bold{y}})) \} \le 0
	\end{equation*} and 	
	\begin{equation*}
		\min_{\bold{y} \in \mathbb B^{|Z|}}\mathcal{FR}(\bold{y})-\mathcal{FR}(\bar {\bold{y}}) \le \min\{0, ( \gamma ( \{(i,j)\},\emptyset)-\mathcal{FR}(\bar {\bold{y}})) \} \le 0
	\end{equation*}	
	hold for any $\bar{\bold{y}} \in \mathcal{Y} \cap \mathbb B^{|Z|}$. This completes the proof.
\end{proof}
Since the $\gamma$ function can be solved in polynomial time based on Observation \eqref{OB}, the inequality \eqref{New_cut}, which contains $|Z|$ functions, can be generated in polynomial time.
\section{Strengthen the Valid Inequality}

The lifting technique of integer programming is a generalization method used to strengthen valid inequalities \cite{AlperLifting2004, BorosLift1975, BorosLift1978, Lifting2000, Lifting1977}. In this section, we explore a method that integrates the general lifting concept with the observations in \ref{OB} to strengthen the inequality \eqref{New_cut}. Specifically, given an  incumbent $\bar{\bold{y}}$, we remain the coefficients $\sum_{(i,j)\in \bar Y}  \min\{0, (\gamma (\emptyset, \{(i,j)\})-\mathcal{FR}(\bar {\bold{y}})) \}$ of the inequality \eqref{New_cut} and perform up-lifting procedures for the associated set $Z \setminus \bar Y$. Specifically, given an incumbent $\bar{\bold{y}}$, we retain the coefficients $\sum_{(i,j)\in \bar Y} \min\{0, (\gamma (\emptyset, {(i,j)})-\mathcal{FR}(\bar{\bold{y}})) \}$ of inequality \eqref{New_cut} and perform up-lifting procedures for the other coefficients associated with the set $Z \setminus \bar Y$. Let $(i,j)^1,(i,j)^2,\dots,(i,j)^{|Z \setminus \bar Y|}$ be an ordering of the edges in $Z \setminus  \bar Y$. In this ordering, given $\bar{\bold{y}} \in \mathcal{Y} \cap \mathbb B^{|Z|}$, we sequentially solve $|Z \setminus  \bar Y|$ lifting problems (see, e.g.,  the section II.2 of \cite{NW88} and the section 2 of \cite{Lifting2000}) to get a new class of valid inequalities
\begin{equation*}
		\theta \ge \mathcal{FR}(\bar {\bold{y}}) + \sum_{(i,j)\in \bar Y}  \min\{0, (\gamma (\emptyset, \{(i,j)\})-\mathcal{FR}(\bar {\bold{y}})) \} (1-y_{(i,j)}) +\sum_{k=1}^{|Z \setminus  \bar Y|} \pi_{(i,j)^k}(\bar {\bold{y}})  y_{(i,j)^k}, 
\end{equation*}
where function $\pi_{(i,j)^r}(\bar {\bold{y}})$ for $r=1,\dots,|Z \setminus  \bar Y|$ denotes the $r$-th exact lifting problem defined by 
 \begin{subequations}\label{lift_P}
	\begin{align} 
		\pi_{(i,j)^r}(\bar {\bold{y}}) = -\mathcal{FR}(\bar {\bold{y}})+ \min~~& \mathcal{FR}({\bold{y}}) -  \sum_{(i,j)\in \bar Y}  \min\{0, (\gamma (\emptyset, \{(i,j)\})-\mathcal{FR}(\bar {\bold{y}})) \} (1-y_{(i,j)}) - \sum_{k=1}^{r-1} \pi_{(i,j)^k}(\bar {\bold{y}})  y_{(i,j)^k} \nonumber \\
		\text{s.t.}~~
		& y_{(i,j)^{r}} = 1   \label{lift_c1}\\ 
		& y_{(i,j)^{k}} = 0, \qquad k = r+1, \dots, |Z \setminus  \bar Y| \label{lift_c2}\\ 
		& \bold{y} \in \mathcal{Y} \cap \mathbb B^{|Z|}. \label{lift_c3}
	\end{align} 
\end{subequations}
Note that the different classes of Problem \eqref{SSP_E} are still embedded within the lifting problem \eqref{lift_P}. This implies that directly solving the exact lifting problem \eqref{lift_P} to obtain a new valid inequality remains challenging. Here, we attempt to derive a modified version of the exact lifting problem \eqref{lift_P}, which allows us to obtain a valid inequality more easily.
\begin{theorem}\label{valid:lift}
	Given an incumbent $\bar{\bold{y}} \in \mathcal{Y} \cap \mathbb B^{|Z|}$, the inequality 
\begin{equation}\label{Lift_valid}
	\theta \ge \mathcal{FR}(\bar {\bold{y}}) + \sum_{(i,j)\in \bar Y}  \min\{0, (\gamma (\emptyset, \{(i,j)\})-\mathcal{FR}(\bar {\bold{y}})) \} (1-y_{(i,j)}) +\sum_{k=1}^{|Z \setminus \bar Y|} \hat \pi_{(i,j)^k}(\bar {\bold{y}})  y_{(i,j)^k}
\end{equation}
	is valid for the RMP \eqref{RMP}, where a modified lifting problem  is 
\begin{subequations}\label{lift_P2}
	\begin{align} 
\hat \pi_{(i,j)^r}(\bar {\bold{y}}) & = \min \Bigl\{ 0, -\mathcal{FR}(\bar {\bold{y}})+ \min \{\mathcal{FR}({\bold{y}}): \eqref{lift_c1},\eqref{lift_c2},\bold{y} \in  \mathbb B^{|Z|}\}  \Bigl\} \label{lift_P2_1}\\
		& = \min \Bigl\{ 0, -\mathcal{FR}(\bar {\bold{y}})+ \gamma \big( \{(i,j)^r\}, \{ (i,j)^{r+1},\dots,(i,j)^{|Z \setminus \bar Y|} \} \big)  \Bigl\} \label{lift_P2_2}
	\end{align} 
\end{subequations}	
for $r = 1,\dots,|Z \setminus \bar Y|$.
\end{theorem}    
\begin{proof}	
	Given $\bar{\bold{y}} \in \mathcal{Y} \cap \mathbb B^{|Z|}$, let $\alpha_{(i,j)}$ be a positive constant that $\alpha_{(i,j)} = -\min\{0, (\gamma (\emptyset, \{(i,j)\})-\mathcal{FR}(\bar {\bold{y}})) \} \ge 0$ for all $(i,j) \in \bar Y$. We aim to find an  inequality 
\begin{equation}\label{Lift_valid2}
	\theta \ge \mathcal{FR}(\bar {\bold{y}}) + \sum_{(i,j)\in \bar Y}  \alpha_{(i,j)} (1-y_{(i,j)}) +\sum_{k=1}^{|Z \setminus \bar Y|}  \beta_{(i,j)^k}  y_{(i,j)^k},
\end{equation}	 
	where $\beta_{(i,j)^k} \in \mathbb R$ is a valid coefficient for a variable $y_{(i,j)^k}$ for all $k=1$ to $|Z \setminus \bar Y|$. We first slightly modify the original lifting problem \eqref{lift_P} into another lifting problem referred to as 
	 \begin{equation} \label{lift_P3}
		\tilde \pi_{(i,j)^r}(\bar {\bold{y}}) =\min \Bigl\{0, -\mathcal{FR}(\bar {\bold{y}})+ \min\{ \mathcal{FR}({\bold{y}}) +  \sum_{(i,j)\in \bar Y} \alpha_{(i,j)} (1-y_{(i,j)}) - \sum_{k=1}^{r-1} \pi_{(i,j)^k}(\bar {\bold{y}})  y_{(i,j)^k}: \eqref{lift_c1}-\eqref{lift_c3}\}  \Bigl\}.
	\end{equation}
	We prove the validity of \eqref{lift_P3} for determing the inequality \eqref{Lift_valid2} using the fundamental concepts of lifting and mathematical induction. For the edge $(i,j)^1$ of an arbitrary ordering, the exact lift problem for the valid coefficient of $y_{(i,j)^1}$ in the inequality \eqref{Lift_valid2} is
\begin{align*}
	& -\mathcal{FR}(\bar {\bold{y}})+ \min\{ \mathcal{FR}({\bold{y}}) +  \sum_{(i,j)\in \bar Y} \alpha_{(i,j)} (1-y_{(i,j)}): \eqref{lift_c1}-\eqref{lift_c3}\} \\ 
	&\ge \min \Bigl\{0, -\mathcal{FR}(\bar {\bold{y}})+ \min\{ \mathcal{FR}({\bold{y}}) +  \sum_{(i,j)\in \bar Y} \alpha_{(i,j)} (1-y_{(i,j)}): \eqref{lift_c1}-\eqref{lift_c3}\}  \Bigl\} =\tilde \pi_{(i,j)^1}(\bar {\bold{y}}),
\end{align*} 
	which provides a valid $\beta_{(i,j)^1} = \tilde \pi_{(i,j)^1}(\bar {\bold{y}})$ for the inequality \eqref{Lift_valid2}.		
	Suppose that we have the valid $\tilde \pi_{(i,j)^k}(\bar {\bold{y}})=\beta_{(i,j)^k}$ for $k=1$ to $r-1$. Given that the coefficients $\beta_{(i,j)^1}=\tilde \pi_{(i,j)^1}(\bar {\bold{y}})$ to $\beta_{(i,j)^{r-1}}=\tilde \pi_{(i,j)^{r-1}}(\bar {\bold{y}})$ in the inequality \eqref{Lift_valid2}​ are valid, the exact lifting problem for determining the valid coefficient of $y_{(i,j)^r}$ in the inequality \eqref{Lift_valid2} is
\begin{align*}
	& -\mathcal{FR}(\bar {\bold{y}})+ \min\{ \mathcal{FR}({\bold{y}}) +  \sum_{(i,j)\in \bar Y} \alpha_{(i,j)} (1-y_{(i,j)}) - \sum_{k=1}^{r-1} \tilde \pi_{(i,j)^k}(\bar {\bold{y}})\}  y_{(i,j)^k}: \eqref{lift_c1}-\eqref{lift_c3}\} \\
	&  \ge \min \Bigl\{0, -\mathcal{FR}(\bar {\bold{y}})+ \min\{ \mathcal{FR}({\bold{y}}) +  \sum_{(i,j)\in \bar Y} \alpha_{(i,j)} (1-y_{(i,j)}) - \sum_{k=1}^{r-1} \tilde \pi_{(i,j)^k}(\bar {\bold{y}})  y_{(i,j)^k}: \eqref{lift_c1}-\eqref{lift_c3}\}  \Bigl\}=\tilde \pi_{(i,j)^r}(\bar {\bold{y}}),
\end{align*} 
which provides a valid $\beta_{(i,j)^r} = \tilde \pi_{(i,j)^r}(\bar {\bold{y}})$ for the inequality \eqref{Lift_valid2}. Since  the value $r \in \{1,\dots,|Z \setminus \bar Y|\}$ and the ordering $(i,j)^1,(i,j)^2,\dots,(i,j)^{|Z \setminus \bar Y|}$ are arbitrary, this confirms the validity of \eqref{lift_P3} for determing the inequality \eqref{Lift_valid2}. Compared to \eqref{lift_P3}, the modified lifting problem \eqref{lift_P2} removes the $\mathcal Y$ of $\eqref{lift_c3}$ and the terms $\sum_{(i,j)\in \bar Y} \alpha_{(i,j)} (1-y_{(i,j)}) - \sum_{k=1}^{r-1} \tilde \pi_{(i,j)^k}(\bar {\bold{y}})  y_{(i,j)^k}$ of the objective function. Since $\bold{y} \in  \mathbb B^{|Z|}$ is the relaxation of $\bold{y} \in  \mathcal {Y} \cap \mathbb B^{|Z|}$ shown in \eqref{lift_c3}, the term $\alpha_{(i,j)} \ge 0$ for all $(i,j) \in \bar Y$, and the term $\tilde \pi_{(i,j)^k}(\bar {\bold{y}}) \le 0$ for all $k=1$ to $|Z \setminus \bar Y|$, we conclude that $\hat \pi_{(i,j)^k}(\bar {\bold{y}})$ is the relaxation of $\tilde \pi_{(i,j)^k}(\bar {\bold{y}})$ and $\hat \pi_{(i,j)^k}(\bar {\bold{y}}) \le \tilde \pi_{(i,j)^k}(\bar {\bold{y}})$ for all $k=1$ to $|Z \setminus \bar Y|$. That is, for a feasible point $(\hat {\bold{y}}, \hat \theta)$, we have
	\begin{align}
	\hat \theta & \ge \mathcal{FR}(\bar {\bold{y}}) + \sum_{(i,j)\in \bar Y}  \alpha_{(i,j)} (1- \hat y_{(i,j)}) +\sum_{k=1}^{|Z \setminus \bar Y|}   \tilde \pi_{(i,j)^k}(\bar {\bold{y}}) \hat y_{(i,j)^k}\nonumber\\
	& \ge \mathcal{FR}(\bar {\bold{y}}) + \sum_{(i,j)\in \bar Y}  \alpha_{(i,j)} (1-\hat y_{(i,j)}) +\sum_{k=1}^{|Z \setminus \bar Y|}   \hat \pi_{(i,j)^k}(\bar {\bold{y}}) \hat y_{(i,j)^k} .\nonumber
\end{align} 
This completes the proof.
\end{proof}
\begin{algorithm}\label{Algo1}
	\SetAlgoLined
	Input: An incumbent $\bar Y$ and an arbitrary ordering $(i,j)^1,(i,j)^2,\dots,(i,j)^{|Z \setminus \bar Y|}$ from $Z \setminus \bar Y$; \\
	Get the value $\mathcal{FR}(\bar {\bold{y}})$;\\
	\For{$\forall (i,j) \in \bar Y$}
	{	
		Solve $\min\{0, (\gamma (\emptyset, \{(i,j)\})-\mathcal{FR}(\bar {\bold{y}})) \}$;\\
	}	
	\For{$r = 1$ to $|Z \setminus \bar Y|$}
	{	
		Solve $\hat \pi_{(i,j)^r}(\bar {\bold{y}}) =  \min \Bigl\{ 0, -\mathcal{FR}(\bar {\bold{y}})+ \gamma \big( \{(i,j)^r\}, \{ (i,j)^{r+1},\dots,(i,j)^{|Z \setminus \bar Y|} \} \big)  \Bigl\}$;\\
	}	
	Output: An inequality \eqref{Lift_valid}.
	\caption{The Separation Problem of \eqref{Lift_valid}}
	\label{alg:1}
\end{algorithm}
We provide Algorithm \ref{alg:1} to solve the separation problem of the inequality \eqref{Lift_valid}. Based on Observation \ref{OB}, Corollary \ref{Algo:poly} demonstrates that Algorithm \ref{alg:1} is a polynomial-time method for the inequality \eqref{Lift_valid}.
\begin{corollary}\label{Algo:poly}
	Algorithm \ref{alg:1} is a polynomial-time method for the inequality \eqref{Lift_valid}.
\end{corollary}
\begin{proof}
	From Observation \eqref{OB}, the function $\gamma (\bar S, \bar N)$ in \eqref{gamma_f} can be solved in polynomial time for two exclusive subsets $\bar S \subseteq Z$ and $\bar N \subseteq Z$. Algorithm \ref{alg:1} solves the $\gamma$ function at most $|Z|$ times within the loops. This completes the proof.
\end{proof}
Although inequality \eqref{Lift_valid} is derived from the modified lifting problem \eqref{lift_P2}, it reveals an intriguing property compared to the inequality \eqref{New_cut}. Specifically, for any ordering $(i,j)^1,(i,j)^2,\dots,(i,j)^{|Z \setminus \bar Y|}$, the inequality \eqref{Lift_valid} is stronger than the inequality \eqref{New_cut} and the $L$-shaped inequality \eqref{L_cut}.
\begin{proposition}\label{prop:strong2}
	Given an incumbent $\bar{\bold{y}} \in \mathcal{Y} \cap \mathbb B^{|Z|}$, the inequality \eqref{Lift_valid} is stronger than the inequality \eqref{New_cut}.
\end{proposition}    
\begin{proof}
	Given $\bar{\bold{y}} \in \mathcal{Y} \cap \mathbb B^{|Z|}$, we set up an arbitrary ordering $(i,j)^1,(i,j)^2,\dots,(i,j)^{|Z \setminus \bar Y|}$ for all elements in $Z \setminus  \bar Y$. Note that given an arbitrary $(i,j)^r \in Z \setminus \bar Y$, we have $\min \{\mathcal{FR}({\bold{y}}): \eqref{lift_c1},\eqref{lift_c2},\bold{y} \in  \mathbb B^{|Z|}\} =  \gamma ( \{(i,j)^r\}, \{ (i,j)^{r+1},\dots,(i,j)^{|Z \setminus \bar Y|} \} )    \ge \gamma ( \{(i,j)^r\}, \emptyset )  $. Therefore, from the equation \eqref{lift_P2_2}, the relation
	\begin{equation*}
		\hat \pi_{(i,j)^r}(\bar {\bold{y}}) =  \{0,-\mathcal{FR}(\bar {\bold{y}})+\min \{\mathcal{FR}({\bold{y}}): \eqref{lift_c1},\eqref{lift_c2},\bold{y} \in  \mathbb B^{|Z|}\}\} \ge \min\{0, \gamma ( \{(i,j)^r\}, \emptyset)-\mathcal{FR}(\bar {\bold{y}})\}
	\end{equation*} 	
	holds for an arbitrary $(i,j)^r \in Z \setminus \bar Y$. This relation provides the result since
	\begin{align*}
		\theta & \ge \mathcal{FR}(\bar {\bold{y}})  
		\sum_{(i,j)\in \bar Y} + \min\{0, (\gamma (\emptyset, \{(i,j)\})-\mathcal{FR}(\bar {\bold{y}})) \} (1-y_{(i,j)}) +\sum_{k=1}^{|Z \setminus \bar Y|} \hat \pi_{(i,j)^k}(\bar {\bold{y}})  y_{(i,j)^k}\\
		&\ge \mathcal{FR}(\bar {\bold{y}}) + \sum_{(i,j)\in \bar Y}  \min\{0, (\gamma (\emptyset, \{(i,j)\})-\mathcal{FR}(\bar {\bold{y}})) \} (1-y_{(i,j)})\nonumber +\sum_{k=1}^{|Z \setminus \bar Y|} \min\{0, (\gamma ( \{(i,j)^k\}, \emptyset)-\mathcal{FR}(\bar {\bold{y}})) \} y_{(i,j)^k}.
	\end{align*}
\end{proof}

\begin{corollary}
	Given an incumbent $\bar{\bold{y}} \in \mathcal{Y} \cap \mathbb B^{|Z|}$, the inequality \eqref{Lift_valid} is stronger than the $L$-shaped inequality \eqref{L_cut}.
\end{corollary}
\begin{proof}
	The result is based on the conclusions of Propositions \ref{prop:strong} and \ref{prop:strong2}.
\end{proof}
Finally, the inequalities \eqref{New_cut} and \eqref{Lift_valid} can generate corresponding integer programming formulations for the set $\mathcal{C}$ of the RMP \eqref{RMP}, allowing us to apply a general cutting plane method in integer programming to solve Problem \eqref{SSP_E} optimally.
\begin{proposition}\label{prop:converge}
	The inequality \eqref{New_cut} or \eqref{Lift_valid} finitely converges to an optimal solution in an exact cutting plane method for Problem \eqref{SSP_E}.
\end{proposition}
\begin{proof}
Theorems \ref{prop:New_valid} and \ref{valid:lift} provide formal proofs of the validity of inequalities \eqref{New_cut} and \eqref{Lift_valid}. Problem \eqref{SSP_E} is equivalent to the integer programming formulations
	\begin{align}
		\min \{\theta: \theta \ge  \mathcal{SSP}( {\bold{\bar y}}) & + \sum_{(i,j)\in \bar Y}  \min\{0, (\gamma (\emptyset, \{(i,j)\})-\mathcal{FR}(\bar {\bold{y}})) \} (1-y_{(i,j)})\nonumber \\ \nonumber 
		& +\sum_{(i,j)\in Z\setminus \bar Y} \min\{0, (\gamma ( \{(i,j)\}, \emptyset)-\mathcal{FR}(\bar {\bold{y}})) \} y_{(i,j)} \; \forall \bar Y \subseteq Z, \bold{y} \in \mathcal{Y}\cap \mathbb B^{|Z|}, \theta \in \mathbb R_+\}  \nonumber
	\end{align} 
	and
	\begin{align}
		\min \{\theta: \theta \ge  \mathcal{SSP}( {\bold{\bar y}}) &+  \sum_{(i,j)\in \bar Y}  \min\{0, (\gamma (\emptyset, \{(i,j)\})-\mathcal{FR}(\bar {\bold{y}})) \} (1-y_{(i,j)})\nonumber \\ \nonumber  
& +\sum_{k=1}^{|Z \setminus \bar Y|} \hat \pi_{(i,j)^k}(\bar {\bold{y}})  y_{(i,j)^k} \; \forall \bar Y \subseteq Z, \bold{y} \in \mathcal{Y}\cap \mathbb B^{|Z|}, \theta \in \mathbb R_+\}.  \nonumber
	\end{align} 
	Since the incumbent solutions $\bold{\bar y} \in \mathcal{Y}\cap \mathbb B^{|Z|}$ are finite, this concludes the proof.
\end{proof}

In conclusion, without Observation  \ref{OB}, we cannot even guarantee that the $L$-shaped inequality \eqref{L_cut} can be generated in polynomial time. Based on the observations, we introduce valid inequalities for Problem \eqref{SSP_E}.  We strengthen these inequalities step by step and demonstrate that these valid inequalities can be generated in polynomial time. These results may offer some potential insights for the exact method on Problem \eqref{SSP_E}. For future research, exploring stronger or facet-defining conditions from other forms of lifting techniques is a suitable direction.

\section*{Acknowledgments}
This work is supported, in part by, NSTC Taiwan 111-2221-E-A49-079 and 111-2221-E-A49-126-MY2.

	\bibliographystyle{abbrv}
	\bibliography{general}

\end{document}